\documentclass[12pt]{article}
\usepackage{amssymb,amsmath}
\usepackage{graphicx}
\usepackage{amsfonts}
\usepackage{float}
\usepackage{flafter}
\usepackage{color}

\labelwidth\leftmargini\advance\labelwidth-\labelsep

\def\R{\relax\ifmmode I\!\!R\else$I\!\!R$\fi}

\def\Z{\relax\ifmmode Z\!\!\!Z\else$Z\!\!\!Z$\fi}

\def\C{\relax\ifmmode C\!\!\!\!I\else$C\!\!\!\!I$\fi}

\def\K{\relax\ifmmode I\!\!K\else$I\!\!K$\fi}

\def\N{\relax\ifmmode I\!\!N\else$I\!\!N$\fi}

\newcounter{defcounter}[section]
\newenvironment{definition}%
{\vspace{0.1cm}\begin{sloppypar}\noindent\stepcounter{defcounter}{\bfseries
Definition
      \thesection.\thedefcounter}}%
{\end{sloppypar}\vspace{0.1cm}}

\newtheorem{theorem}{Theorem}[section]

\newtheorem{proposition}{Proposition}[section]

\newcommand{\proof}{{\bf Proof.} }

\newcommand{\qed}{\hfill $\square$}

\begin{document}
\thispagestyle{empty}
\begin{center}
{\Large {\bf Non-uniform expansions of real numbers}}
\end{center}
\begin{center}J\"org Neunh\"auserer\\
joerg.neunhaeuserer@web.de
\end{center}
\begin{center}
\begin{abstract}
We introduce and study non-uniform expansions of real numbers, given by two non-integer bases.
\end{abstract}
\end{center}
\section{Introduction} Expansions of reals numbers in non-integer bases are studied since
the pioneering works of R\'{e}nyi in the end of the 1950s and Parry in the 1960s, see \cite{[RE],[PA1],[PA2]}. In these works especially the greedy algorithm  that determines the digits of such expansions and the relationship of these expansions to symbolic dynamics is addressed. In the 1990s a group of Hungarian mathematics led by Paul Erd\H{o}s revived this file of research, see \cite{[E1],[E2],[E3]}. Beside other results they proved that each $x\in(0,1/(1-q))$ has a continuum of expansions of the form $\sum_{i=1}^{\infty} q^{-n_{i}}$ if $1<q<G$, where $G$ is the golden ratio. In the sequel Sidorov \cite{[SI]} used ergodic theoretical methods to prove that for all $q\in(1,2)$ almost all  $x\in(0,1/(1-q))$ have such an expansion. Moreover, Glendinning and Sidorov \cite{[GS]} proved that there always exist (at least countably many) reals having a unique expansion if $q>G$. Nowadays especially dimensional theoretical aspects of expansions of reals numbers in non-integer bases are studied, see for instance \cite{[KK],[KK2],[BA]}.\\
In this paper we introduce non-uniform expansions of real numbers, which may be viewed as expansions with respect to two non-integer bases. As far as we know such expansions were not studied yet, although they constitute a natural generalisation. The rest of the paper is organized as follows:\\
In the next section we give two descriptions of non-uniform expansions of real numbers. In the following we introduce a greedy, a lazy and intermittent algorithms that give the digits of these expansions. In section four we prove a theorem on the existence of a continuum of non-uniform expansions of real numbers, which is similar to the results in the uniform case we mentioned above. In the last section we characterise real numbers which have a unique non-uniform expansion and prove a theorem on the cardinality of the set of such numbers.
\section{The expansions}
Let $\Sigma=\{0,1\}^{\mathbb{N}}$ be the set of sequences of zeros and ones. Equipped with the metric
\[ d((s_{i}),(t_{i}))=\sum_{i=1}^{\infty}|s_{i}-t_{i}|2^{-i}\]
$\Sigma$ is a compact, perfect and totally disconnected space. For $s=(s_{i})\in\Sigma$ and $n\in\mathbb{N}$ let $0_{n}(s)$ be the number of zeros and $1_{n}(s)$ be the number of ones in the sequence $(s_{1},\dots,s_{n})$.\\
Fix $\beta_{0},\beta_{1}\in (1/2,1)$ with $\beta_{0}\ge \beta_{1}$ and let $I=I_{\beta_{1}}=[0,\beta_{1}/(1-\beta_{1})]$.\footnote{In the literature the uniform case $\beta_{0}=\beta_{1}$ has been studied. Usually the reciprocal of $\beta_{0}$ resp. $\beta_{1}$ is denoted by $\beta$.}
We consider the map $\pi_{\beta_{0},\beta_{1}}:\Sigma \to I$ given by
\[ \pi_{\beta_{0},\beta_{1}}(s)=\pi_{\beta_{0},\beta_{1}}((s_{i}))=\sum_{i=1}^{\infty}s_{i}\beta_{0}^{0_{i}(s)}\beta_{1}^{1_{i}(s)}.\]
Several times we will use another description of this map, which we now describe. Let $T_{0},T_{1}:I\to I$ be the contractions given by
\[ T_{0}(x)=\beta_{0} x\qquad\mbox{ and }\qquad T_{1}(x)=\beta_{1}x+\beta_{1}.\]
By induction we have
\[ T_{s_{1}}\circ\dots \circ T_{s_{n}}(x)=\beta_{0}^{0_{n}(s)}\beta_{1}^{1_{n}(s)}x+\sum_{i=1}^{n}s_{i}\beta_{0}^{0_{i}(s)}\beta_{1}^{1_{i}(s)}.\]
Hence for all $s\in\Sigma$ and all $x\in I$
\[ \pi_{\beta_{0},\beta_{1}}(s)=\pi_{\beta_{0},\beta_{1}}((s_{i}))=\lim_{n\to\infty}T_{s_{1}}\circ\dots \circ T_{s_{n}}(x).\]
\begin{definition}
We call a sequence $s\in\Sigma$ with $\pi_{\beta_{0},\beta_{1}}(s)=x$ a $(\beta_{0},\beta_{1})$-expansion of $x\in I$.
\end{definition}~\\
The following proposition guarantees the existence of $(\beta_{0},\beta_{1})$-expansions.
\begin{proposition}
The map $\pi_{\beta_{0},\beta_{1}}$ is continuous and surjective.
\end{proposition}
\proof
If  $d((s_{i}),(i_{i}))<2^{-u}$ we have $s_{i}=t_{i}$ for $i=1,\dots ,u$, which implies
\[ |\pi_{\beta_{0},\beta_{1}}(s_{i})-\pi_{\beta_{0},\beta_{1}}(t_{i})|<\beta_{0}^{u}\beta_{1}/(1-\beta_{1}).\]
Hence $\pi_{\beta_{0},\beta_{1}}$ is continuous.
Note that
\[ T_{0}(I)\cup T_{1}(I)
=[0,\beta_{0}\beta_{1}/(1-\beta_{1})]\cup[\beta_{1},\beta_{1}/(1-\beta_{1})]=I\]
since $\beta_{0}+\beta_{1}\ge 1$. This implies
\[ \bigcup_{(s_{i})\in\{0,1\}^n} T_{s_{1}}\circ\dots \circ T_{s_{n}}(I)=I.\]
for every $n\ge1$. Hence for each $x\in I$ there is sequence $s\in \Sigma$ such that
\[ x\in  T_{s_{1}}\circ\dots \circ T_{s_{n}}(I),\]
but this implies $\pi_{\beta_{0},\beta_{1}}(s)=x$. Therefore $\pi_{\beta_{0},\beta_{1}}$ is surjective
\qed~\\~\\
In the next section we describe an algorithm which determines one $(\beta_{0},\beta_{1})$-expansion of $x\in I$.
\section{The greedy and the lazy algorithm}
Using the notations of last section we define a a map $G:I\to I$
\[ G(x)=\left\{\begin{array}{ll} T_{0}^{-1}(x), & x\not\in T_{1}(I) \\
         T_{1}^{-1}(x), & x \in T_{1}(I)  \end{array}\right.
  \]
\[=\left\{\begin{array}{ll} \beta_{0}^{-1}x, & x\in [0,\beta_{0}) \\
         \beta_{1}^{-1}x-1, & x \in [\beta_{0},\beta_{1}/(1-\beta_{1})],  \end{array}\right.. \]
For $x\in I$ we define the greedy expansion $g=g(x)=(g_{i})\in\Sigma$ with respect to $(\beta_{1},\beta_{2})$ by
\[ g_{i}=\lfloor\beta_{1}^{-1}G^{i-1}(x)\rfloor,\]
where $\lfloor a\rfloor$ is the greatest interner not greater than $a$. We have
\begin{proposition}
For all $x\in I$ the greedy expansion $g(x)\in\Sigma$ with respect to $(\beta_{0},\beta_{1})$ is an $(\beta_{0},\beta_{1})$-expansion of $x$, that is $\pi_{\beta_{0},\beta_{1}}(g(x))=x$.
\end{proposition}
\proof
If $g_{i}=0$ we have $G^{i-1}(x)<\beta_{1}$,  which implies $G^{i-1}(x)\not\in T_{1}(I)$ and $G^{i-1}(x)\in T_{g_{i}}(I)$. If $g_{i}=1$ we have $G^{i-1}(x)\ge\beta_{1}$,  which implies $G^{i-1}(x)\in T_{1}(I)$ hence again $G^{i-1}(x)\in T_{g_{i}}(I)$. By the definition of $G$ we conclude
\[ x\in  T_{g_{1}}\circ\dots \circ T_{g_{i}}(I)\]
for all $i\in\mathbb{N}$,  but this implies $\pi_{\beta_{0},\beta_{1}}(g(x))=x$.
\qed~\\~\\
\begin{figure}
\vspace{0pt}\hspace{50pt}\scalebox{0.5}{\includegraphics{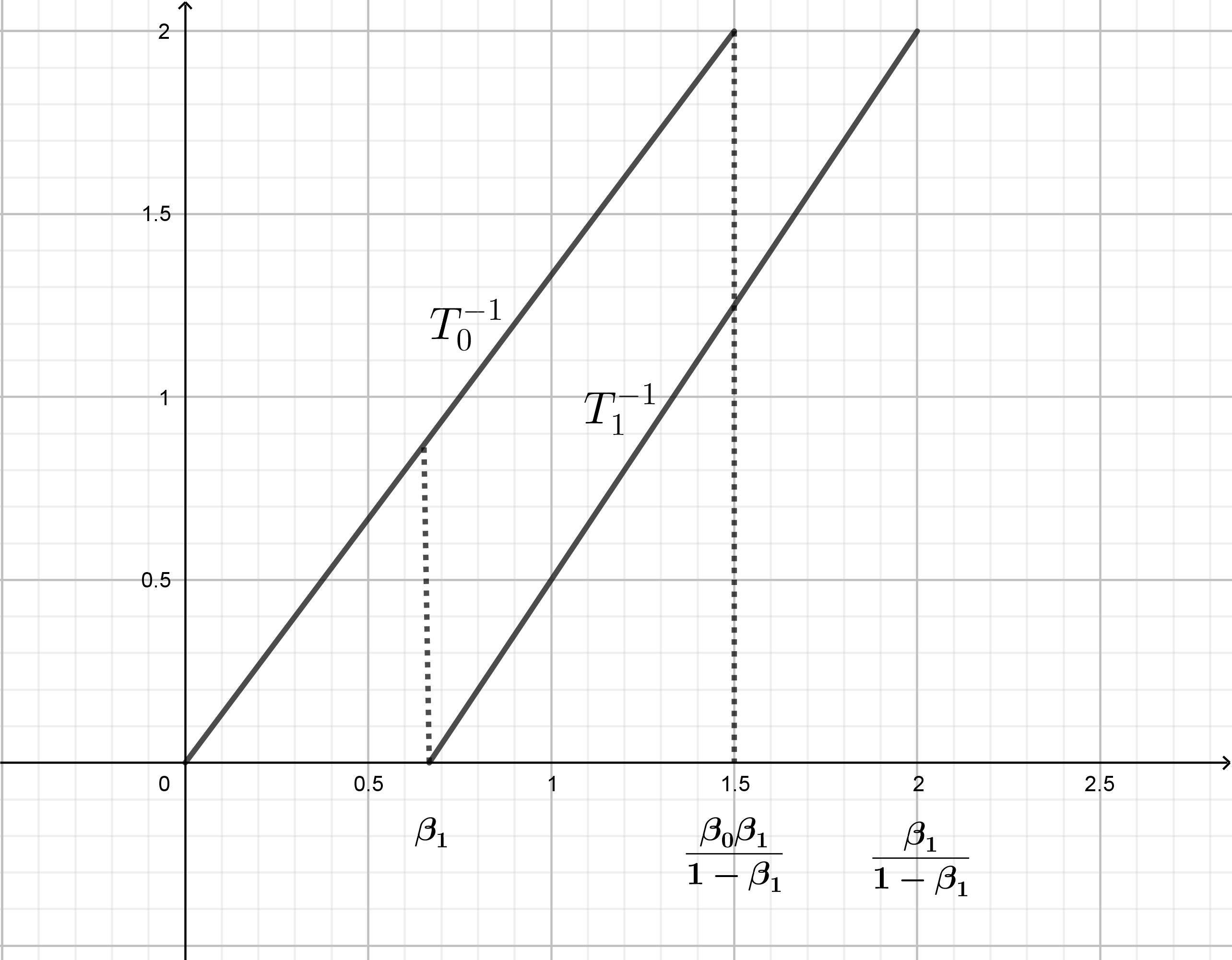}}
\caption{The maps $T_{0}^{-1}$ and $T_{1}^{-1}$ in the case $\beta_{0}=3/4,\beta_{1}=2/3$}
\end{figure}
To define the lazy expansion let $L:I\to I$ be given by
\[ L(x)=\left\{\begin{array}{ll} T_{0}^{-1}(x), & x\in T_{0}(I) \\
         T_{1}^{-1}(x), & x \not\in T_{0}(I)  \end{array}\right.
  \]
\[=\left\{\begin{array}{ll} \beta_{0}^{-1}x, & x \in [0,\beta_{0}\beta_{1}/(1-\beta_{1})] \\
         \beta_{1}^{-1}x-1, & x \in (\beta_{0}\beta_{1}/(1-\beta_{1}),\beta_{1}/(1-\beta_{1})],  \end{array}\right. ~.\]
For $x\in I$ the lazy expansion $l=l(x)=(l_{i})\in\Sigma$ with respect to $(\beta_{1},\beta_{2})$ is given by
\[ l_{i}=\lceil(1-\beta_{1})(\beta_{0}\beta_{1})^{-1} L^{i-1}(x)\rceil-1,\]
where $\lceil a\rceil$ is the smallest integer not smaller than $a$.
We have
\begin{proposition}
For all $x\in I$ the lazy expansion $l(x)\in\Sigma$ with respect to $(\beta_{1},\beta_{2})$ is an $(\beta_{0},\beta_{1})$-expansion of $x$, that is $\pi_{\beta_{0},\beta_{1}}(l(x))=x$.
\end{proposition}
\proof
If $l_{i}=0$ we have $L^{i-1}(x)\le\beta_{0}\beta_{1}/(1-\beta_{1})$,  which implies $L^{i-1}(x)\in T_{0}(I)$ hence $L^{i-1}(x)\in T_{l_{i}}(I)$. If $l_{i}=1$ we have $L^{i-1}(x)>\beta_{0}\beta_{1}/(1-\beta_{1})$,  which implies $L^{i-1}(x)\not\in T_{0}(I)$ hence again $L^{i-1}(x)\in T_{l_{i}}(I)$. By the definition of $L$ we conclude
\[ x\in  T_{l_{1}}\circ\dots \circ T_{l_{i}}(I)\]
for all $i\in\mathbb{N}$,  which implies $\pi_{\beta_{0},\beta_{1}}(l(x))=x$.
\qed~\\~\\
For $\alpha\in (\beta_{1},\beta_{0}\beta_{1}/(1-\beta_{1})$ we may also consider intermediate expansions $m(x)=(m_{i})$ with respect to $(\beta_{1},\beta_{2})$ given by
\[ m_{i}=\lfloor\alpha^{-1}M_{\alpha}^{i-1}(x)\rfloor,\]
where
\[M_{\alpha}(x)=\left\{\begin{array}{ll} \beta_{0}^{-1}x, & x\in [0,\alpha) \\
         \beta_{1}^{-1}x-1, & x \in [\alpha,\beta_{1}/(1-\beta_{1})]\end{array}\right..\]
Again these are  $(\beta_{0},\beta_{1})$-expansions of $x\in I$.
\section{A continuum of expansions}
It is natural to ask how many $(\beta_{0},\beta_{1})$-expansions a real number in $I$ has. It turns out that usually there is a continuum of such expansions:
\begin{theorem}
Let $\beta_{0},\beta_{1}\in (1/2,1)$ with $\beta_{0}\ge \beta_{1}$. We have:\\
(1) Almost all $x\in I$ have  a continuum of $(\beta_{0},\beta_{1}$)-expansions.\\
(2) If $\beta_{1}^2+\beta_{0}> 1$ all $x\in I\backslash\{0,\beta_{1}/(1-\beta_{1})\}$ have a continuum of $(\beta_{0},\beta_{1})$-expansions.\\
\end{theorem}
\proof We first prove (2). Let $J=(0,\beta_{1}/(1-\beta_{1}))$ and
\[\Lambda_{0}=T_{0}(J)\cap T_{1}(J)=(\beta_{1},\beta_{0}\frac{\beta_{1}}{1-\beta_{1}}).\]
We recursively define $\Lambda_{n+1}=T_{0}(\Lambda_{n})\cup\Lambda\cup T_{1}(\Lambda_{n})$ and
prove by induction:
\[ \Lambda_{n}=(\beta_{0}^{n}\beta_{1},(\beta_{1}^{n}(\beta_{0}-1)+1)\frac{\beta_{1}}{1-\beta_{1}}).\]
We have
\[T_{0}(\Lambda_{n})\cup\Lambda\cup T_{0}(\Lambda_{n})=(\beta_{0}^{n+1}\beta_{1},(\beta_{1}^{n}\beta_{0}(\beta_{0}-1)+1)\frac{\beta_{1}}{1-\beta_{1}})\cup  (\beta_{0}^{n}\beta_{1},(\beta_{1}^{n}(\beta_{0}-1)+1)\frac{\beta_{1}}{1-\beta_{1}})\]\[\cup                                           (\beta_{0}^{n}\beta_{1}^2+\beta_{1},(\beta_{1}^{n+1}(\beta_{0}-1)+1)\frac{\beta_{1}}{1-\beta_{1}})=\Lambda_{n+1}.\]
In the last equation we use $\beta_{1}^2+\beta_{0}>1$ and $\beta_{0}^2+\beta_{1}> 1$, which is true since $\beta_{0}\ge \beta_{1}$. Note that $\bigcup_{n\ge 0} \Lambda_{n}=J$. Hence for every $x\in J$ there is a $k\ge 0$ and a sequence $(s_{1},\dots,s_{k})\in\{0,1\}^{k}$ such that
\[ x=T_{s_{1}}\circ\dots\circ T_{s_{k}}\circ T_{0}(x_{0})\mbox{ and }x=T_{s_{1}}\circ\dots\circ T_{s_{k}}\circ T_{1}(x_{1}),\]
where $x_{0},x_{1}\in J$ and $x_{0}\not=x_{1}$. Hence we obtain two expansions of $x$ that differ in the $k+1$-digit. Applying the result to $x_{0}(x)$ and $x_{1}(x)$ we obtain four expansions of $x$. Here we use that $x_{0}(x)$ and $x_{1}(x)$ are not at the boundary of $J$. Repeating this procedure $\aleph_{0}$ times we see that there are $2^{\aleph_{0}}$ expansions of $x$.~\\
Now we prove (1). Let $G:I\to I$ be the map associated with the greedy expansion from section 3. $G$ is a piecewise linear expanding interval map and such maps are known to have an ergodic measure, which is equivalent to the Lebesgue measure, see \cite{[DE]} and \cite{[KH]}. By Poincare recurrence theorem for almost all $x\in I$ there is a $k\ge 0$ such that $G^{k}(x)\in \Lambda_{0}$. Hence for almost all $x\in J$ there is a $k\ge 0$ and a sequence $(s_{1},\dots,s_{k})\in\{0,1\}^{k}$ such that
\[ x=T_{s_{1}}\circ \dots\circ T_{s_{k}}\circ T_{0}(x_{0})\mbox{ and }x=T_{s_{1}}\circ\dots\circ T_{s_{k}}\circ T_{1}(x_{1}),\]
where $x_{0},x_{1}\in J$ and $x_{0}\not=x_{1}$. For almost all $x$ both numbers $x_{1}(x),x_{2}(x)$ have two different $(\beta_{0},\beta_{1})$-expansion hence almost all $x$ have four different expansions. We use here that the intersection of two sets of full measure has full measure. Repeating this procedure $\aleph_{0}$ times we obtain $2^{\aleph_{0}}$ expansions for almost all $x\in I$, using the fact a countable intersection of sets of full measure has full measure.
\qed~~\\~\\
Obviously the $(\beta_{0},\beta_{1})$-expansion of $0$ and $\beta_{1}/(1-\beta_{1})$ is unique. Our theorem leaves the question open if there are numbers $x$ in the interior of $I$ that have a unique $(\beta_{0},\beta_{1})$-expansion. We will address this question in the following section.
\section{Unique expansions}
We consider the shift map $\sigma:\{0,1\}^{\mathbb{N}}\to \{0,1\}^{\mathbb{N}}$ given by $\sigma((s_{k}))=(s_{k+1})$. Using this map we may characterise numbers which have a unique $(\beta_{0},\beta_{1})$-expansion as follows:
\begin{proposition}
The $(\beta_{0},\beta_{1})$-expansion $(s_{i})$ of $x$ is unique if and only if
\[ \pi_{\beta_{0},\beta_{1}}(\sigma^{k}(s_{i}))\in [0,\beta_{1})\cup(\beta_{0}\beta_{1}/(1-\beta_{1}),\beta_{1}/(1-\beta_{1})]\]
for all $k\ge 0$.
\end{proposition}
\proof $\pi_{\beta_{0},\beta_{1}}((s_{i}))=\pi_{\beta_{0},\beta_{1}}((t_{i}))$ with $(s_{i})\not=(t_{i})$ if and only if there exists a smallest $k\ge 0$ such that $s_{k}\not=t_{k}$ and $\pi_{\beta_{0},\beta_{1}}(\sigma^{k}(s_{i}))=\pi_{\beta_{0},\beta_{1}}(\sigma^{k}(t_{i}))$. But this is equivalent to $ \pi_{\beta_{0},\beta_{1}}(\sigma^{k}(s_{i}))\in T_{0}(I)\cap T_{1}(I)=[\beta_{1},\beta_{0}\beta_{1}/(1-\beta_{1})]$. The proposition follows by contraposition.\qed~\\~\\
Using this characterisation of points with unique expansion we are able to prove:
\begin{theorem}
Let $\beta_{0},\beta_{1}\in (1/2,1)$ and $\beta_{0}\ge \beta_{1}$. \\
(1) If $\beta_{0}(1+\beta_{1})<1$ there are at least countable many $x\in I$, which have a unique $(\beta_{0},\beta_{1})$-expansion.\\
(2) If $\beta_{0}(1+2\beta_{1}-\beta_{0}\beta_{1})<1$ there are uncountable many $x\in I$, which have a unique $(\beta_{0},\beta_{1})$-expansion. Moreover the set of these $x$ has positive Hausdorff dimension.
\end{theorem}
\proof First we prove (1). Consider the periodic sequence $p=(010101\dots)$. Since $\beta_{0}(1+\beta_{1})<1$ we have
\[ \pi_{\beta_{0},\beta_{1}}(p)=\beta_{0}\beta_{1}/(1-\beta_{0}\beta_{1})<\beta_{1}. \]
Note that $\beta_{0}(1+\beta_{1})<1$ implies $\beta_{1}(1+\beta_{0})<1$ since $\beta_{0}\ge \beta_{1}$. Hence we have $\beta_{0}-\beta_{0}^2\beta_{1}<1-\beta_{1}$ and thus
\[ \pi_{\beta_{0},\beta_{1}}(\sigma(p))=\beta_{1}/(1-\beta_{0}\beta_{1})=\beta_{0}\beta_{1}/(\beta_{0}-\beta_{0}^2\beta_{1})>\beta_{0}\beta_{1}/(1-\beta_{1}).\]
By proposition $x=\pi_{\beta_{0},\beta_{1}}(p)$ has a unique $(\beta_{0},\beta_{1})$-expansion. Obviously the same is true for all $x$ of the form
$x=\pi_{\beta_{0},\beta_{1}}((0\dots 0101010\dots))$ and there countable many of such $x$.\newline
Now we prove (2). Let $V=\{01,10\}^{\mathbb{N}}$ and
\[ U=\bigcup_{k=0}^{\infty}\sigma^{k}(V)=V\cup(\{0\}\times V)\cup \cup(\{1\}\times V).\]
We prove that $\pi_{\beta_{0},\beta_{1}}(U)\subseteq  [0,\beta_{1})\cup(\beta_{0}\beta_{1}/(1-\beta_{1}),\beta_{1}/(1-\beta_{1})]$.
The sequence $s\in U$ with $s_{1}=0$ that has the largest projection under $\pi_{\beta_{0},\beta_{1}}$ obviously is $s=(011010101\dots)$. We have
\[  \pi_{\beta_{0},\beta_{1}}(s)=\beta_{1}\frac{\beta_{0}+\beta_{0}\beta_{1}-\beta_{0}^2\beta_{1}}{1-\beta_{0}\beta_{1}}<\beta_{1}\]                                                            by our assumption. The sequence $s\in U$ with $s_{1}=1$ that has the smallest projection under $\pi_{\beta_{0},\beta_{1}}$ obviously is $s=(1001010101\dots)$. We have
\[  \pi_{\beta_{0},\beta_{1}}(s)=\beta_{1}+\frac{(\beta_{0}\beta_{1})^2}{1-\beta_{0}\beta_{1}}>\frac{\beta_{0}\beta_{1}}{1-\beta_{1}}.\]
The inequality here is equivalent to $\beta_{1}(1+2\beta_{0}-\beta_{0}\beta_{1})<1$ which is true since we assume $\beta_{0}\ge \beta_{1}$. It remains to show that the Hausdorff dimension of $A:=\pi_{\beta_{0},\beta_{1}}(V)$ is positive. Consider the maps
\[ F(x)=T_{0}\circ T_{1}(x)=\beta_{0}\beta_{1}x+\beta_{0}\beta_{1}\]
and
\[ G(x)=T_{1}\circ T_{0}(x)=\beta_{0}\beta_{1}x+\beta_{1}\]
 and let $J=(0,\beta_{1}/(1-\beta_{0}\beta_{1})$. We have $F(J)\subseteq J$ and $G(J)\subseteq J$. Moreover
 \[ F(J)\cap G(J)=(0,\beta_{0}\beta_{1}^2/(1-\beta_{0}\beta_{1}))\cap (\beta_{1},\beta_{1}/(1-\beta_{0}\beta_{1}))=\emptyset \]
 by our assumptions on $\beta_{0}$ and $\beta_{1}$. In the language of fractal geometry this means that $(F,G)$ induce an iterated function system fulfilling the open set condition, see \cite{[FA]}. The attractor of this iterated function system is $A$ since $A=F(A)\cup G(A)$ and the classical formula for self-similar fractals gives
 \[ \dim_{H}A=\frac{-\log(2)}{\log(\beta_{0}\beta_{1})}>0.\]
 \qed

\end{document}